\documentclass[12pt]{amsart}

\usepackage{latexsym}
\usepackage{amsmath}
\usepackage{amsfonts}
\usepackage{amssymb}
\usepackage{graphicx}

\newtheorem{theorem}{Theorem}[section]

\theoremstyle{definition}
\newtheorem{definition}[theorem]{Definition}
\newtheorem{example}[theorem]{Example}

\theoremstyle{remark}
\newtheorem{remark}[theorem]{Remark}

\numberwithin{equation}{section}

\begin{document}

\begin{center}
{\large \bf AROUND FINITE-DIMENSIONALITY IN FUNCTIONAL ANALYSIS-\\ A PERSONAL PERSPECTIVE}\\

{\bf by \vskip .5em M.A. Sofi} \\

Department of Mathematics, Kashmir University, Srinagar-190006,
INDIA email: aminsofi@rediffmail.com

\end{center}

{\bf Abstract\footnote{Mathematics Subject Classification (2000)
\\
Primary 46A11, 46C15; Secondary 46B20, 47A68;
}:}  As objects of study in functional
analysis, Hilbert spaces stand out as special objects of study as do
nuclear spaces (ala Grothendieck) in view of a rich geometrical
structure they possess as Banach and Frechet spaces, respectively.
On the other hand, there is the class of Banach spaces including
certain function spaces and sequence spaces which are distinguished
by a poor geometrical structure and are subsumed under the class of
so-called Hilbert-Schmidt spaces. It turns out that these three
classes of spaces are mutually disjoint in the sense that they
intersect precisely in finite dimensional spaces. However, it is
remarkable that despite this mutually exclusive character, there is
an underlying commonality of approach to these disparate classes of
objects in that they crop up in certain situations involving a
single phenomenon-the phenomenon of finite dimensionality-which, by
definition, is a generic term for those properties of Banach spaces
which hold good in finite dimensional spaces but fail in infinite
dimension.

\section{Introduction}

\noindent A major effort in functional analysis is devoted to
examining the extent to which a given property enjoyed by finite
dimensional (Banach) spaces can be extended to an infinite
dimensional setting. Thus, given a property (P) valid in each finite
dimensional Banach space, the desired extension to the infinite
dimensional setting admits one of the following three (mutually
exclusive)
possibilities:\\
(a). (P) holds good in all infinite dimensional Banach spaces.\\
(b). (P) holds good in some infinite dimensional Banach spaces.\\
(c). (P) fails in each infinite dimensional Banach space.\\
Examples of (P) verifying (a) include, for instance, the validity of
the inverse function theorem that one learns in an advanced course
on calculus whereas a typical example of property (P) verifying (b)
is provided by reflexivity or say, the Radon-Nikodym property of a
Banach space. Our main focus will be on property (P) falling under
(c)-henceforth to be designated as finite-dimensional properties-
which will be seen to exhibit remarkable richness when considered in
an infinite dimensional setting. This latter setting will naturally
necessitate the consideration of an infinite dimensional framework
beyond the world of Banach spaces in which these properties are
sought to be salvaged, which in our case will be provided by the
class of Frechet spaces. Interestingly, this leads to the
identification of nuclear (Frechet) spaces as natural objects in
which- at least some of the most important -finite-dimensional
properties are valid in an infinite dimensional setting outside the
framework of Banach spaces. On the other hand, as an important
subclass of Banach spaces, there is the class of so-called
Hilbert-Schmidt spaces which include certain important function
spaces like $ C(K),L_{1} (\mu),L_{\infty} (\mu)$ or sequence spaces
like $c_{0}, \ell_{\infty}, \ell_{1}$ besides many more. These
spaces are characterized by the property that a bounded linear
operator acting between Hilbert spaces and factoring over such a
space is already a Hilbert Schmidt operator. In particular, no
infinite dimensional Banach space can simultaneously be a Hilbert
space as well as a Hilbert -Schmidt space. In a similar vein, a
Hilbert space can never be nuclear unless it is finite dimensional.
It turns out that, along with nuclear spaces, these three (disjoint)
classes of spaces which are distinguished by properties that are
mutually exclusive arise as different manifestations of this single
phenomenon of finite dimensionality. The present article is devoted
to a discussion of certain interesting aspects of this phenomenon
and to point out the ways in which it leads to the consideration of
nuclear spaces on the one hand and of Hilbert spaces and
Hilbert-Schmidt spaces on the other.

\section{Notation}
\noindent Throughout this paper, we shall let X, Y, Z denote Banach
spaces, unless otherwise stated. We shall use the symbol $X^{*}$ for
the dual of X whereas $B_{X}$ shall be used for the closed unit ball
of X:

$$ B_{X} = \{x\in X: \| x\|\leq1\}.$$
We shall also make use of the following notation:
\begin{center}
$L(X,Y)$, Banach space of all bounded linear operators from $X$ into
$Y$.
\end{center}
\begin{center}
 $P(\mathbb{N})$, set of all permutations on $\mathbb{N}$.
\end{center}
\begin{center} $c_{0}(X)=\{( x_{n} ) \subset X:  x_{n}\rightarrow 0\}
$
\end{center}
$$\ell[X ] = \bigg\{( x_{n} ) \subset X:\sum _{n=1}^{\infty }{|<x_{n},f>|}^{p}<\infty,\forall f\in X^{*}\bigg\}$$
$$c_{0}(X)=\{( x_{n} ) \subset X:  x_{n}\rightarrow 0\} $$
$$\ell_{p}[X ] = \bigg\{( x_{n} ) \subset X:\sum _{n=1}^{\infty }{||<x_{n}>||}^{p}<\infty,\bigg\}.$$
Clearly, $\ell_{p}\{X\}  \subset\ell_{p}[X ], \forall  1\leq
p<\infty $. Further, the indicated inclusion is continuous when
these spaces are equipped with natural norms defined by:
$$\varepsilon\bigg((x_{n})\bigg ) = \sup_{f\in B_{X^{*}}}  \sum _{n=1}^{\infty }{|<x_{n},f>|}^p, (x_{n})\in \ell_p[X].$$
$$\pi\bigg((x_{n})\bigg ) = {\bigg(\sum _{n=1}^{\infty }{\|x_{n}\|}^{p}\bigg)}^{1/p}, (x_{n})\in \ell_{p}\{X\}.$$
For a locally convex space $X$, a locally convex topology on each of
these sequence spaces can be defined exactly in the same fashion,
with the norm being replaced by a family of seminorms generating its
topology.

\noindent For later use, let us note the following equivalent
description of these spaces of sequences identifying them as certain
classes of bounded linear maps.
\begin{theorem}
 For $1<p<\infty $ and $q$ where $\frac{1}{p}+\frac{1}{q}=1$ the correspondence\\
$T\rightarrow \{ T(e_p)\} $  provides an isometric isomorphism of\\
(i)$L(\ell_p,X) \,\,onto\,\, \ell_p[X]$\\
(ii)$\Pi_p(\ell_p, X)\,\,onto\,\,\ell_p\{X\}$.
\end{theorem}
\noindent Here $\Pi_p$ denotes the class of p-summing maps defined
below (see Definition 3). The proof of (i) is a straightforward
consequence of the fact that $||T||=||T^*|| $ and that $||T^*|| $ is
nothing but $\varepsilon_p\bigg((x_n)\bigg)$.

\noindent An outline of a proof of (ii) is provided in Theorem
3.17(a).
\section{Definitions and Examples}
\noindent Given a property $(P)$, we say that $(P)$ is a {\it
finite-dimensional property} $((FD)$-property, for short) if it
holds good for all  finite dimensional Banach spaces but fails for
each infinite dimensional Banach space.
\begin{example}
(i)Heine-Borel Property (closed bounded subsets of X are compact).\\
(ii) $X^* = X'$ (algebraic dual of X).\\
(iii) $ w^{*}$-convergence = norm convergence in $X^*$ (Josefson-Nissenweig theorem) \\
(iv) Completeness of the weak-topology on X.\\
(v) Equivalence of all norms on X. \\
(vi) Hahn-Banach property: For a given Banach space X, the following
holds: \\
Each bounded linear operator defined on a subspace of X and taking
values in an arbitrary Banach space Y can be extended to a bounded
linear operator on X and each bounded linear operator defined on X
can be extended to a bounded linear operator on any superspace of
X.\\
(vii) Dvoretzky-Hanani Property: Given a null sequence $x_n$ in X,
there exist signs $\epsilon_n = \pm 1, n\geq 1$ such that
$\sum_{n=1}^{\infty }\epsilon_n x_n  $
 converges in X.\\
(viii) McArthur Property:A series $\sum_{n=1}^{\infty } x_n  $ in X
such that $\sum_{n=1}^{\infty } x_{\pi(n)}  $
 converges to the same sum for all permutations $\pi \in P(\mathbb{N})$ for which
 $\sum_{n=1}^{\infty } x_{\pi(n)}$ converges is already
unconditionally convergent.\\
(ix)Riemann-Rearrangement Property(RRP)(also called Dvoretzky-Rogers
property (DRP)):
$$uv(X) = \ell_{1}[X]= \ell_{1}\{X\} = abc(x).$$
Here $ uc(X) = \bigg \{(x_{n}) \subset X:\sum_{n=1}^{\infty }
x_{\pi(n)}$ Converges in X,$ \forall \pi \in P(\mathbb{N})\bigg\} $
and
 $$abc(X) = \bigg\{( x_n)\subset X:\sum_{n=1}^{\infty }||x_n||<\infty \bigg\}.$$

 \noindent Clearly, $abc(X)\subset uc(X)$. More generally, we have the equality:
 $$\ell_p\{X\}=\ell_p[X], 1 \le p<\infty $$

as an important (FD)-property. This is the p-analogue of the famous
Dvoretzky- Rogers theorem referred to above.
\end{example}
 \begin{definition} ( [9], Chap. 2): T $\in L(X,Y)$ is said to be\\
 (i) nuclear (T$\in N(X,Y)$ if there exist
 $\{ \lambda_n\} \subset \ell_1, \{ f_n\} \subset X^* \, and \, \{ f_n\} \subset Y$)
  bounded such that $$T(x) = \sum_{n=1}^{\infty }\lambda_n<x,f_n>y_n, x\in X$$
 The class $N(X,Y)$ then becomes a Banach space when equipped with the nuclear norm defined by
  $$\nu(T) = inf \bigg\{\sum_{n=1}^{\infty }|\lambda_n|\,||f_n||\,||y_n||\bigg\},$$
 where infimum is taken over all representations of T as given above.\\

 (ii) p- ($\it absolutly)$ summing $\bigg(T \in  \Pi_p(X,Y)\bigg) \,
 if \,\forall \{x_n\} \in \ell_p[X],$  it follows that
 $\bigg\{(T(x_n)\bigg\} \in  \ell_p\{X\}.$ By the open mapping theorem,
 this translates into the finitary condition: $\exists \,c > 0$ s.t.
\begin{equation}
  {\bigg(\sum_{i=1}^{n } {||x_i||}^{p}\bigg)}^{1/p}\leq c. \,
 \sup_{f\in B_{X^*}}\bigg\{ {\bigg(\sum_{i=1}^{n }{|<x_i,f>|}^p\bigg)}^{1/p}\bigg\},
\end{equation}
$\forall\, {(x_i)}^{n}_{i=1} \subset X,  n \geq 1 .$
 \noindent The infimum of all such c$'s$ appearing above and denoted by $\pi_{p} (T)$ is called the $p$-{\it summing norm} of T, making $\Pi_p (X,Y)$ into a Banach space.
\end{definition}
\noindent Grothendieck's theorem (see \cite{24}, Theorem 5.12)
states that all
 operators on an $L_{1}$ space and taking values in a Hilbert space are absolutely
  summing whereas those acting on an $L_{\infty}$ space are always 2-summing.
 \begin{remark} (See [22], Chap.II).\\
(i) Nuclear maps are always compact (as the uniform limit of a sequence of finite rank operators)\\
(ii) $N(X,Y)\subset \Pi_p (X,Y)\subset \Pi_q (X,Y), \,\,\,\forall 1\leq p\leq  q.$\\
(iii)$\Pi_p (H_1,H_2)= HS(H_1, H_2)$, class of Hilbert-Schmidt maps
acting between Hilbert
spaces $H_1$ and $H_2$ and $1\leq p < \infty $.\\
(iv) $T \in \Pi_p (X,Y)$ if and only if $TS \in \Pi_p (\ell_q ,Y)$
for each $S \in  L(\ell_q, X), \frac{1}{p}+\frac{1}{q}=1.$\\
 (v) A composite of 2-summing maps (and hence of absolutely summing maps) is always nuclear with
 $\nu (TS) \leq \pi_{2}(T)\pi_{2}(S), $ for all
 $S\in \Pi_{2}(X,Y), T\in \Pi_{2}(Y,Z)$.
\end{remark}
\noindent {\bf Three important features of (FD)-properties:}\\
As stated in the Introduction, the commonality of approach to the
three classes of spaces, namely nuclear spaces, Hilbert spaces and
Hilbert-Schmidt spaces through finitedimensional phenomena is based
on three important features lurking in the shadows of these
phenomena but which manifest themselves only in an infinite
dimensional context. These three features involving a given finite
dimensional property (P) derive from:\\
a. Frechet space analogue of (P).\\
b. Size of the set of objects failing (P).\\
c. Factorisation property of (P).\\

\noindent In view of its impact on the development of functional
analysis via Dvoretzky-Rogers theorem and a host of other variants
of this latter theorem, in what follows we shall mainly discuss the
(RRP) as an important (FD)-property and single it out as the main
example to illustrate the aforementioned phenomena while making an
effort to provide the necessary details in respect of other
(FD)-properties, wherever possible. Before we proceed further, let
us pause to argue why (RRP) is indeed an (FD)-property. To this end,
we recall:\\ Dvoretzky-Rogers Lemma ([9], Lemma 1.3): Given $n \geq
1$ and a 2n-dimensional normed space X, there exist n vectors
$x_{1}, x_{2}\cdots x_{n}\in B_{X}$ with $||x_{i}||\geq 1/2$ for
$1\leq n$ such that
above lemma to assert the following $(*)$ which is a stronger
statement than that $(RRP)$ is an $(FD)$ property. \\ $(*)$ Given an
infinite dimensional Banach space X and
${(\lambda_{i})}_{i=1}^{\infty} \in \ell_{2},$ there exists
${(x_i)}_{i=1}^{\infty}\subset X$  such that
$\sum_{n=1}^{\infty}x_n$ is unconditionally convergent in X and
 $||x_i||= |\lambda_i|, i\geq 1.$

\noindent Indeed, we can choose an increasing sequence
${(n_k)}_{k=1}^{\infty}$ of positive integers such that for m, n ,
one has
$$\sum_{i=m}^{n}{|\lambda_i|}^{2} \leq {2}^{-2k} $$

\noindent Applying the above Dvoretzky-Rogres lemma to each block
$[n_k, n_{k+1})$ yields a sequence ${(y)}_{i=1}^{\infty}\subset
B_{X}$ with $||y_i||\geq 1/2$ for $i \geq 1$ such that for all
choices of scalars $\alpha_i$ and $n_k \leq N\le n_{k+1},$ , we have
$$\| \sum_{i=n_{k}}^{N} \alpha_i y_i \|\leq {\bigg(\sum_{i=n_k}^{N} {|\alpha_i|}^{2}\bigg)}^{1/2}  $$

Setting $x_{i} = \lambda_{i} y_{i}/||y_{i}||$ and taking
$\epsilon_{i} = \pm 1$  for $i\geq 1$ gives

$$\| \sum_{i=n_k}^{N} \epsilon_i x_i\|\leq {\bigg(\sum_{i=n_k}^{N} \frac{{|\lambda_i|}^2}
{{\|y_i\|}^2}\bigg)}^{1/2} \leq 2^{-k+1}, \,\, k\geq 1.$$

\noindent In other words, the series $\sum_{i=n_k}^{N} \epsilon_i
x_i$ has Cauchy partial sums and is, therefore, convergent.
Equivalently, the series $\sum_{i=1}^{\infty} x_n$  is
unconditionally convergent $\|
x_i \| = |\lambda_i|$ for for all $i\geq 1.$\\

\noindent {\bf Comments:} Besides the proof of the Dvoretzky-Rogers
(DR) - theorem given above, there are many more, exploiting
techniques ranging from local theory of Banach spaces to topological
tensor products. Instead, we shall briefly include details of yet
another proof of this important statement which is based on a
theorem of the author proved in [27]. It characterizes the finite
dimensionality of a Banach space X in terms of the equality of
certain operator ideals:
\begin{theorem}
 The following statements for a Banach space X are equivalent: \\
(i) $\Pi_{2} (X, \ell_{2}) = N(X, \ell_{2}).$\\
(ii) dim$X\le \infty$.

\noindent To see how the proof of the (DR)-theorem based on Theorem
2 works, we note that the given condition in (DR)-theorem, namely:
$\ell_1 [X] =\ell_1 \{X\} $ leads to the equality:
$\ell_{2}[X]=\ell_{2}\{X\}$which, by virtue of Theorem 2.1
translates into:L$(\ell_{2},X)=\Pi_{2}(\ell_{2},X)$.Thus there
exists $c>0$ such that
\begin{equation*}
 \pi_{2}(T) \leq c\parallel T \parallel,\forall T \in L(\l_{2},X)
 \end{equation*}
Now a trace duality argument applied to the above equality gives:
$\Pi_{2}(\ell_{2},X)=N(X,\ell_{2})$,
 which is the same thing as X being finite dimensional, by Theorem 3.4.
Indeed, given ܶ$T\in \Pi_{2}(\ell_{2},X)$ and finite dimensional
spaces E and F with $u\in L(E,X),v\in L(\l_{2},F),w\in  L(F,E)$,we
note that $uvw\in L(\ell_{2},X)$, and so, applying (3.1)to $uvw$ and
combining the resulting estimate with Remark 3.3 (v)gives:
\begin{equation*}
 \nu (uwvT)\leq \pi_{2}(uwv)\pi_{2}(T)\leq c\parallel  uvw\parallel \pi_{2}(T)
\end{equation*}
This together with an application of [9], Lemma 6.14 yields:
\begin{equation*}
  trace(wvTu)= trace(uwvT)
\end{equation*}
This gives:
\begin{equation*}
 i(T)\leq c\pi_{2}(T)
\end{equation*}
where i denotes the ‘integral norm’. In other words, T is an
integral operator taking its values in a reflexive space, and so has
to be nuclear. It follows that  $\Pi_{2}(X,\ell_{2})=N(X,\ell_{2})$
and, therefore, X is finite dimensional by Theorem 3.4.
\end{theorem}

\noindent {\bf (a) Frechet-space setting:}\\ When suitably
formulated in the setting of Frechet spaces X, it turns out that in
most of the cases, there exist infinite dimensional Frechet spaces
in which an (FD)- property holds. It also turns out that, at least
in most cases of interest, the class of Frechet spaces in which this
holds coincides with the class of nuclear spaces in the sense of
Grothendieck \cite{11}. To define a nuclear space, we shall assume-
in the interest of technical simplifications- that the topology of
locally convex spaces in question will be generated by a system of
norms as opposed to a family of seminorms.
\begin{definition} A locally convex space(lcs, for short) X is said to
be nuclear if for each continuous norm p on X, there exists a
continuous norm q on X ,$ q \geq  ‫$   p   ‬ such that the
identity map i:(X,q)$ \longrightarrow$(X,p) is nuclear. It is easily
verified that this is equivalent to requiring that for each Banach
space Z, each continuous linear map T: X$\longrightarrow$Z ܼ is
nuclear.
\end{definition}
\begin{example}
The following are well known examples of nuclear spaces;
\begin{enumerate}
\item  All finite dimensional spaces.
\item $\omega$ countable product of the line.
\item H(C), space of entire functions(on the plane).
\item H(D), space of holomorphic functions on the unit disc D.
\item D($\Omega$), space of test functions on an open set $\Omega$ in $\mathbb{R}^{n}$ .
\item D'($\Omega$) , space of distributions on an open set $\Omega$ in  $\mathbb{R}^{n}$.
\end{enumerate}
\end{example}
\begin{remark}
Banach+Nuclear = Finite Dimensional.
\end{remark}
\begin{example} The Frechet- space analogue of the following
(FD)-properties are valid exactly when the Frechet space in question
is nuclear.
\end{example}
\noindent (i) (RRP): Unconditionally convergent series in X are absolutely convergent.\\
(ii) Levy-Steinitz Property: Given a convergent series
$\sum_{n=1}^{\infty}$ in X, then DS($\sum_{n=1}^{\infty}$) is an
afine space. more precisely,
\begin{equation*}
 DS \bigg(\sum_{n=1}^{\infty}x_{n} \bigg)=\Gamma \bigg(\sum_{n=1}^{\infty}x_{n} \bigg)+\sum_{n=1}^{\infty}x_{n}
\end{equation*} where
$\Gamma \bigg(\sum_{n=1}^{\infty} x_{n} \bigg)= = \{ x\in X;
f(x)=0,\forall f \in X^{*} s.t. \\ \sum_{n=1}^{\infty} |< x_{n},f>|
\textless \infty \}$ and $DS \bigg(\sum_{n=1}^{\infty}x_{n} \bigg)$,
the domain of sums of the series $\sum_{n=1}^{\infty}x_{n}$is
defined by:
\begin{equation*}
 DS \bigg(\sum_{n=1}^{\infty}x_{n} \bigg)=\bigg\{  x\in X; \exists    \pi
  P(\mathbb{N})s.t \sum_{n=1}^{\infty}x_{\pi(n)}=x \bigg\}
\end{equation*}
(iii)  $ D(C (K, X), Y) = \Pi_{1} ((C (K, X), Y).$ Here,$D (C (K,
X), Y)$ stands for the class of dominated operators: $T߳\in D (C
(K, X), Y)$ if there exists a Borel measure $\mu$ on K such that

\begin{equation*}
 \| Tf \|\leq\int_{k}\| f(t)\| d\mu (t), \forall f\in C(K,X)
\end{equation*}
(iv) Bochner Property(BP): Positive definite functions f on X arise
as Fourier transforms of regular Borel measures on $X^{*}$ :

\begin{equation*}
 f(x)= \int_{X^{*}} e^{-i<x,x^{*}>} d\mu, x\in X
\end{equation*}
(v) Weakly closed subgroups of a Frechet space X coincide with
closed subgroups of X.\\
(vi)  Equality involving vector measures:\\
(a)$M(X)=M_{b\nu}(X).$\\
Here, $M(X ) \quad and \quad M_{b\nu}(X)$ denote the spaces of
X-valued measures:
$$   M(X) =  \bigg\{\mu: A  \rightarrow X:\mu \bigg(\bigcup_{n=1}^{\infty} A_{n}\bigg)
= \sum_{n=1}^{\infty}\mu (A_{n}), A_{m} \cap A_{n}= = \phi, \forall
m \neq n\bigg\} $$
$$   M_{bv}(X) =  \bigg\{\mu \in M(X):
 \sup_{(A_n)\subset \mathcal{A}}
 \sum_{n=1}^{\infty} \|\mu (A_n)\|< \infty, A_m \cap A_n = \phi, \forall m \neq n\bigg\} $$
(b) $c_0(X) \subset R_{bv} (X)$\\
where the symbols involved have the following meanings:\\
$ R(X) = \{ (x_n) \subset X:\exists \mu \in M(X),(x_n) \subset rg(\mu)\} $\\
$ R_{bv}(X) = \{ (x_n) \subset X:\exists \mu \in M_{bv}(X),(x_n) \subset rg(\mu)\} $\\
$rg(\mu) = \{ \mu(A): A \in A\} $\\

\noindent {\bf Comments:} Let us begin by saying that each of the
above properties which are valid in nuclear Frechet spaces provide
yet further evidence that as opposed to Hilbert spaces, nuclear
spaces are more suited to be looked upon as infinite dimensional
analogues of finite dimensional spaces. Further, the proof for the
equivalence of each of the above properties to nuclearity is
different in each case, drawing upon techniques from different areas
of analysis, depending upon the nature of the property. We shall,
however, settle for a sketch of proof of the equivalence of (1) with
nuclearity, using by now the standard techniques from the theory of
operator ideals proposed by Pietsch as opposed to the complicated
approach via  tensor products
which was orginally given by Grothendieck in his thesis [11].\\

\noindent Proof of (i): Let us say that operator ideals $A$ and $B$
are equivalent if some power of  $A$ is contained in $B$ and some
suitable power of $A$ is contained in ऋ. The definition of nuclear
spaces given earlier shows that the class of locally convex spaces
determined by an operator ideal that is equivalent to the ideal of
nuclear operators in the sense just described coincides with the
class of nuclear spaces. Now, (ii) and (iv) of Remark 2(a) shows
that the ideals N and ߎଵ are equivalent. Combining this fact with
the definition of an absolutely summing map completes the proof.
Indeed, assume that each unconditionally convergent series in the
Frechet space X is absolutely convergent. This yields that the
inclusion $i: \ell_1\{X\}\rightarrow \ell_1[X]$ is a well-defined
bijective continuous (linear) map. Further, the completeness of X
yields that each of these spaces is complete w.r.t. the metrisable
locally convex topologies defined in (I). Thus, the inverse mapping
theorem applies to yield that the inverse map $i = i^-:
\ell_1[X]\rightarrow \ell_1\{X\}$ is continuous. In other words, for
each continuou norm p on X, there exists a continuous norm q on X,
$q\geq p$ and $c >0$ such that
$$ \sum_{i=1}^{n}p(x_{i})\leq c. \sup _ {f\in B_{X^{*}_{q}}}
\bigg(\sum_{i=1}^{n}|< x_{n} f>|\bigg),\forall (x_{i})^{n}_{i=1}
\subset X, n\geq 1.$$

\noindent By 3.2(ii), this means that the identity map $i:
(X,q)\rightarrow (X,p)$‬is absolutely summing and this completes
the argument. Converse is a strightforward consequence of 3.3(i)\\

\noindent (ii) The first complete proof of this statement was given
by W.Banasczyk \cite{5} in 1990 which makes use of certain
combinatorial lemmas now knows as the ``Rearragement Lemma'' and the
``Lemma on Rounding off co-efficients'' involving a finite set of
vectors in a (metrisable) nuclear space. A further strengthening of
this statement valid in all complete (DF) - spaces was proposed by
Bonet and Defant [7] a little later. A unified treatment of the
(LS)-theorem which subsumes the finite dimensional as well as the
nuclear analogue besides certain instances of its validity in a
Banach space setting was given by author in \cite{30}.\\

\noindent (iii) The result is folklore for X, a 1-dimensional space
and generalizes easily to the case of X being finite dimensional. C.
Swartz \cite{32} proved the converse of this latter statement, i.e.,
the equality $D(C(K,X),Y)  $ for each Banach space Y is an
(FD)-property for X. The Frechet space analogue of the equality was
established by M.Nakamura ሾ19ሿ who showed the equivalence of
this equality to the nuclearity of X. Here, the necessity part is a
direct consequence of the definition of a nuclear space whereas
sufficiency can be proved by using (i) above which consists in
showing that each unconditionally convergent series in X is
absolutely convergent. Thus, let $\sum_{n=1}^{\infty} x_n$ be an
unconditionally convergent series in X and fix t in K. Consider the
map $T: C(K,X)\rightarrow X$‬ defined by $T(f) = f(t).$ ‬Letting
$\{p_n;n\geq 1\}$ denote a generating family of (semi)norms for the
topology of X and $\delta_t$ the Dirac measure at t, we see that for
each $n\geq 1, p_n(T(f))=<p_n(f(),\delta_t>$ which shows that T is
1-dominated, and hence absolutely summing by the given hypothesis.
Let $ g\in C(K)$ be the function: $g(s) =1,$ for all s in $K$. By
identifying the dual of ‫$C(K,X)$ with $X^*$ - valued regular
(c.a.) Borel measures on K, it follows that the sequence $\{g( )x_n;
n\geq 1\}$   is weakly absolutely summable in $C(K, X) $, and so
$\{T(g( ) x_n; n\geq 1\}= \{x_n; n\geq 1\}$  is weakly absolutely
summable in $X,$  i.e., the series $\sum_{n=1}^{\infty} x_n$
absolutely convergent.\\

\noindent (iv) The necessity part of Bochner’s theorem in the
setting of lcs was proved by Minlos \cite{17} who showed its
validity in metrisable nuclear spaces. The converse that the
validity of Bochner’s theorem in a metrisable lcs X implies
nuclearity of X is a deep result of D. Muschtari \cite{18}.\\

\noindent (v) It is a well- known theorem of S.Mazur that the weak
closure of a convex set in a normed space coincides with its
nom-closure. In particular, a (linear) subspace of a normed space
$X$ is closed if and if it is weakly closed. However, this
equivalences does not carry over to subgroups of an infinite
dimensional normed space X. In fact, the existence of such groups is
important from the viewpoint of unitary representations of
topological groups as it leads to easy examples of topological
groups which may not even admit weakly continuous unitary
representations. The existence of closed subgroups which are not
weakly closed was proved proved by S. J. Sidney \cite{25} in the
setting of infinite dimensional Banach spaces admitting a separable
infinite dimensional quotient whereas the general case for an
infinite dimensional normed space was settled by Banaszczyk
\cite{2}. Also, the validity of this property in nuclear
(metrizable) spaces is again due to Banaszczyk \cite{3} and the
equivalence of this property with the nuclearity of a metrisable lcs
is due to M. Banasczyk and W. Banasczyk \cite{4}. As a consequence,
it follows that, as in the case of locally compact groups, nuclear
Frechet spaces admit unitary representations on a Hilbert space
which is even faithful if the topology of the space is given by a
sequence of norms.\\

\noindent (vi) In the theory of vector measures, It is known that
each vector measure taking its values in a finite dimensional space
is of bounded variation. A simple application of Dvoretzky-Rogers
theorem yields that there are vector measure taking values in an
infinite dimensional Banach space which are not of bounded
variation. Same is true of (b) which says that the property
involving the containment of null sequences from a Banach space X
inside the range of X-valued measures is an (FD)-property. Duchon
\cite{10} showed that nuclear spaces are the only Frechet spaces X
for which each vector measure taking values in X is of bounded
variation whereas the equivalence of (b) with the nuclearity of X
was proved by Bonet and Madrigal \cite{8}. As a strengthening of the
latter property, it was shown by the author \cite{28} that this
equivalence remains valid for the smaller sequence space
$\ell_p\{X\}$ in place of $c_0(X)$\\

\noindent The above discussion serves to bring home the view that,
as opposed to Hilbert spaces, nuclear spaces provide the most
convenient infinite dimensional framework for the validity of
certain important (FD)-properties which, by their very definition,
fail in each infinite dimensional Banach space.\\

\noindent {\bf(b).Size of the set of objects failing (FD)}\\
1. Given an (FD)- property (P), it turns out that for a given
infinite- dimensional
 Banach space X, the set of objects in X failing (P) is usually ‘very big’: it could be topologically big(dense)\\
algebraically big(contains an infinite-dimensional space)\\
big in the sense of category(non-meagre),\\
big in the sense of functional analysis (contains an
infinite-dimensional closed subspace).\\

\noindent 2. {\underline{Examples}}:\\
(i). For an infinite-dimensional Banach space $X$, the difference
set $uc(X)/ abc(X)$ contains a c-dimensional
subspace.\\
(ii) $M(X)/ M_{bv}(X)$ is non-meagre.\\
(iii) $D(C (K), X) / \Pi (C (K), X)$ contains an
infinite-dimensional space.\\
(iv) $M([0,1], X) /B([0,1], X)$ contains an infinite-dimensional space.\\
Here, M and B stand, respectively, for the class of McShane and
Bochner
integrable functions on $[0, 1]$ taking values in X.\\
(v) $DS(\sum_{n=1}^{\infty}$ the {\it domain of sums} in the Levy-Steinitz theorem is far from being convex.\\

\noindent {\bf Comments:} As in (a) above, let us see how to prove
(i) which measures the extent of failure of the Riemann
Rearrangement Theorem in an infinite dimensional Banach space.\\

\noindent {\underline{Proof of (i)}}: We begin by considering an
uncountable almost disjoint family $\{A_{\alpha}\}_{\alpha \in
\Lambda}$ infinite subsets of  $\mathbb{N}: \, A_{\alpha}\cap
A_{\beta}$ is a finite set for $\alpha \neq \beta $. An easy way to
produce such a family is by writing $\Lambda = [0,1]$ and letting
$\{r_n\}^{\infty}_{k=1}$ denote the rationals in $[0,1]$. For
$\alpha \in \Lambda$, choose a subsequence $\{
r_{n_{k}}\}^{\infty}_{k=1}$ of $\{
r_{n_{k}}\}^{\infty}_{n=1}$ such that $r_{n_{k}}\begin{matrix} k\\
\longrightarrow \end{matrix} \alpha $ and define $A_{\alpha}= \{
n_k: n \geq 1\}$. By construction, the family
$\{A_{\alpha}\}_{\alpha} \in \Lambda $ has the desired properties.
By Dvoretzky-Rogers theorem, we can choose a series
$\sum_{n=1}^{\infty} x_n$
 in $X$ which is unconditionally convergent but which does not converge absolutely. For every $\alpha \in \Lambda$
define a sequence $x^{\alpha}={(x^{\alpha}_i)}^{\infty}_{i=1} $ in
$X$ where $x^{\alpha}=x_n$ if $i = nth$ element of $A_{\alpha}$ and
$x^{\alpha}=0$ otherwise. Clearly, the series $\sum_{n=1}^{\infty}
x_i^{\alpha}$
 is unconditionally convergent.
However, it is not absolutely convergent as it has a subseries which
is not absolutely convergent. By virtue of almost disjointness of
$A'_{\alpha}s$, it follows that the family $\{x^{\alpha}: \alpha \in
\Lambda\}$ is linearly independent. Thus $E= span
\{x^{\alpha}:\alpha\in \Lambda\} $ is c- dimensional
($c=$cardinality of continuum). Noting that each element of $E$ is
unconditionally summable, the proof will be completed by showing
that each element of $E$ is absolutely non-summable. To this end,
let $\{\alpha_1, \alpha_2, \cdots \alpha_n\}\subset \Lambda$ and
$\lambda_1,\lambda_2,\cdots, \lambda_n $ be (non-zero) scalars.
Then, $z=\lambda_1x^{\alpha_1},\lambda_2x^{\alpha_2},
\cdots\lambda_nx^{\alpha_n} $ is not absolutely summable. Indeed, we
can choose an infinite set $A\subset A_{\alpha_1}$ such that
$A_{\alpha_1}- A$ is finite, $A\cap
({U}^{n}_{i=2}A_{\alpha_i})=\phi$ and, therefore, $\sum_{i\in
A}z_{i} = \sum_{i\in A}\lambda_{1} x^{i}_{\alpha_{1}}$ which
contains all but finitely many terms of a non-absolutely convergent
series. It follows that $\sum_{i=1}^{\infty}z_i$ is a non-absolutely
convergent series.

\noindent Property (ii) was proved by R. Anantharaman and K.M.Garg
\cite{1}. However, it is not known if this set (together with 0)
also contains an infinite dimensional space. In his report, one of
the referees has suggested a strategy drawing on some recent
techniques in the area of ‘spaceability’ that would yield that
the set under reference contains a large infinite dimensional space!
Also (iii) is a recent result of F.J.G.Pacheco and Puglisi \cite{20}
whereas the proof of (iv) is part of joint work of the author with
F.J.G.Pacheco \cite{21}.\\
(v) The first example demonstrating the failure of Levy-Steinitz
theorem in infinite dimensional Hilbert space was given by
Marcinkiewicz which can be used, via Dvoretzky’s spherical
sections theorem, to show that such examples can also be constructed
in each infinite dimensional Banach space. The strategy involved in
Marcinkiewicz’s example consists in showing that the constructed
series has a nonconvex domain of sums and, therefore, such a series
fails the Levy-Steinitz theorem. Although in the case under
reference, it is not clear what the “size of the set of objects
failing the Levy-Steinitz” should mean, it is possible to quantify
the extent of failure of this theorem in terms of the degree of
non-convexity of the domain of sums. In fact, it is possible to
produce, in each infinite dimensional Banach space, counterexamples
to the Levy-Steinitz theorem for which the domain of sums is {\it
highly non-convex} in the sense that it consists precisely of two
distinct points$!$\\

\noindent {\bf(c) . Factorization of (FD)-properties}\\
1. As opposed to Hilbert spaces which possess a rich geometrical
structure, the class of those Banach spaces which lie at the other
end of the spectrum is distinguished by a relatively poor
geometrical structure as is testified by spaces like $c_0,
\ell_{\infty}, \ell_1, C(K), L_{\infty}(\Omega) $ etc. The latter
class of Banach spaces is subsumed under the class of so-called
Hilbert-Schmidt spaces which arise as spaces with a particular
factorisation property. To put this definition into perspective, let
us recall that a trace class(nuclear) operator $T:\ell_{2}
\rightarrow \ell_{2} $  admits a factorization over each
infinitedimensional Banach space $X$: there exist bounded linear
maps $T_1:\ell_{2} \rightarrow X, \,\, T_{2} : X\rightarrow
\ell_{2}$ such that $T= T_{2}T_{1}$. Even more, this statement also
holds for the more general class of Hilbert-Schmidt maps.
Conversely, a Hilbert space operator factoring over each infinite
dimensional Banach space is necessarily of the Hilbert-Schmidt type.
For this, it is enough to choose a test space for X which may be
taken to be any of the spaces listed above:
$c_{0},\ell_{\infty},\ell_{1}, C(K), L_{\infty}(\Omega) $ besides
many more that include the disc algebra Î(î). This motivates the
following definition, due originally to H.Jarchow \cite{13}.

\begin{definition}
 $X$ is said to be a Hilbert-Schmidt space if
each bounded linear operator acting between Hilbert spaces and
factoring over X is already a Hilbert Schmidt map.
\end{definition}

\begin{remark}
Recalling the well-known extension property of
2-summing maps, it is possible to prove the following
characterization of Hilbert-Schmidt spaces:
\end{remark}

\begin{theorem}
 For a Banach space X, the following are equivalent:\\
(i) $X$ is a Hilbert-Schmidt space\\
(ii) $L(X, \ell_{2}) =\Pi_{2}(X, \ell_{2} )$.
\end{theorem}
\noindent Proof: The simple proof of $(i) \Rightarrow (ii)$ is a
consequence of the definition of a Hilbert-Schmidt space combined
with the fact that a map $T\in L(X,Y)$ is p-summing if and only if
$TS$ is p-summing for each $S \in L(\ell_p,X)$. The other
implication follows from Remark $ 3.3(iii)$.

\begin{example}
(i) Hahn-Banach Extension Property: For an extremally disconnected
space K and an arbitrarily given Banach space $X$, every bounded
linear operator on $C(K)$ into $X$ extends to a bounded linear
operator on any superspace containing $C(K)$. (Converse is also
true). Something similar is true for all the spaces listed above:
$c_0,\ell_{\infty},\ell_{1},$\\ $C(K),A(D),  L_{\infty}(\Omega) $: a
bounded linear operator defined on any of these spaces and taking
values in $\ell_{2}$ extends to a bounded linear operator on any
Banach space containing the given space. In fact, we have the
following theorem which is a consequence of Theorem c.3 combined
with the extension property of 2-summing maps. The reverse
implication follows by noting that every Banach space embeds into a
C(K) space and that by Grothendieck’s theorem (see Section II), an
$\ell_{2}$-valued bounded linear map on a C(K)-space is already
2-summing.
\end{example}
\begin{theorem}
For a given Banach space X, TFAE:\\
(i) $X$ is a Hilbert-Schmidt space.\\
(ii) $A$ bounded linear map on X into $\ell_{2}$ extends to a
bounded linear map on any superspace of $X$.
 \end{theorem}

\noindent Let us note a dual counterpart to this result:

\begin{theorem}
 For a Banach space $X$ with $dim> 2$, TFAE:\\
(i) A bounded linear map defined on a subspace of $X$ and taking
values in an
arbitrary Banach space Z extends to a bounded linear map on $X$.\\
(ii) $X$ is a Hilbert space.
\end{theorem}

\noindent The proof of $(ii) \Rightarrow (i)$  follows from the
projection theorem in Hilbert spaces whereas $(i) \Rightarrow (ii)$
is a consequence of the famous Lindenstrauss-Tzafriri theorem
\cite{16} to the effect that projection theorem holds good precisely
when the space in question is Hilbertian.

\begin{remark}
 The (FD)-property involving a multivariate analogue of the
HahnBanach extension property lends itself to a similar
factorisation scheme as has been witnessed in the case of the
classical Hahn-Banach in Theorems 3.13 and 3.14 above. More
precisely, we have the following theorem on the extension of
bilinear forms on a Banach space:
 \end{remark}

\begin{theorem}
 A Banach space $X$ is a Hilbert space if and only if for each
2-dimensional subspace $Y$ of $X$, every continuous bilinear
functional on
$Y\times Y^*$  extends to a continuous bilinear functional on $X\times Y^*$ \\
Theorem(b)[14]. A Banach space $X$  such that for each Banach space
$Z$  containing $X$ isometrically, each continuous bilinear form on
$X\times X$ extends to a continuousbilinear form on $Z\times Z$ is a
Hilbert-Schmidt space.
\end{theorem}
\noindent The fact that a Hilbert space can never be Hilbert-Schmidt
space unless it is finite dimensional motivates the following
problem:\\

\noindent {\bf Factorization Problem:} Given an (FD)-property(P),
whether it is possible to write $(P) = (Q) \wedge (R)$ for some
properties $(Q)$ and $(R)$ such that $X$ verifies $ (Q)$ iff $X$ is
Hilbertian and $X$ verifies $(R$) iff $X$ is Hilbert-Schmidt.

\noindent As noted above, the property (P) defined by the
``Hahn-Banach extension property on subspaces {\underline{and}}
superspaces'', admits such a factorization. In what follows, we
shall come across a number of (FD) - properties which admit such a
factorization,some of which are listed below.\\
\noindent (i) $\ell_{2}\{X\}=\ell_{2}[X] $\\
(ii)Bochner Property(BP)\\
(iii) $c_{0}(X) \subset R_{bv}(X)$\\
(iv) $\Pi_{2}(X, .)= \Pi^{d_{2}}(X, .)$\\
(v)$\Pi_{2}(X, \ell_{2})= N(X,\ell_{2})$\\

\noindent {\bf  Missing Link:} There is a 'missing link' involved in
each of the above properties (P) which when inserted between the
objects involved gives rise to properties (Q) and (R) with $(P) =
(Q) \wedge (R)$. Let us see what this missing link looks
like in case of (i) above.\\

 \noindent(i)The (FD)-property in question
involves the equality of two X-valued sequence spaces for which one
sided inclusion is always true: $\ell_{2}\{X\}\subset \ell_{2}[X], $
regardless of the space $X$. The missing object in the factorization
of the equality in (i) entails the search for an $X$-valued sequence
space
 $\lambda(X)$ for which it always holds that $\ell_{2}\{X\}\subset \lambda(X) \subset \ell_{2}[X]$
 and that strengthening these inclusions to equalities results in X being (isomorphically) a Hilbert
  space and a Hilbert-Schmidt space, respectively. In the case under study, the right candidate for
   $\lambda(X)$ turns out to be the space:
$$\ell_{2}\{X_{0}\}= \bigg\{ (x_{n})\subset X: \sum_{n=1}^{\infty}{\|T(x_{n})\|}^{2} < \infty, \,\,
\forall T\in L(X,\ell_{2})\bigg\}$$ for which it holds that
$\ell_{2}\{X\}\subset\ell_{2}\{X_{0}\}\subset\ell_{2}[X] $. Now, the
desired factorization of the equality in (i) is given in the
following theorem:
\begin{theorem}
For a Banach space X, the following statements are true:\\
(a) $\ell_{2}\{X\}=\ell_{2}\{X_{0}\}$ if and only if $X$ is Hilbertian.\\
(b)$\ell_{2}\{X_{0}\}=\ell_{2}[X]$ if and only if X is a
Hilbert-Schmidt space.
\end{theorem}
\noindent {\bf Proof (a):} The proof is accomplished by noting that
the vector-valued sequence spaces appearing here can be identified
with the space $\Pi_{2}(\ell_{2}$ of 2-summing maps and the space
 $\Pi^{d}_{2}(\ell_{2},X)$ of dual 2-summing maps, respectively. The desired correspondence is provided by
 the map $T \rightarrow \{T(e_{n})\}$ which sets up an isometric isomorphism between\\
(i)$\Pi_{2}(\ell_{2},X) \, \, and \, \, \ell_{2}\{X\} $ and between\\
(ii)$\Pi^{d}_{2}(\ell_{2},X) \, \, and \, \, \ell_{2}\{X_{0}\} $

\noindent To see why it is so, let us begin by remarking that the
proof of (i) as given below can be suitably modified to prove the
p-analogue of (i) as stated in Theorem 2.1. Now a simple consequence
of the Hahn-Banach theorem yields that $x= (x_n)\in \ell_{2}\{X\}$
if and only if there exists $a = (a_{n})\in \ell_{2}$ such that
$|<x^{*},x_{n}>| \leq a_{n},\,\, \forall n \geq 1$ and $\forall
x^{*} \in B_{X^{*}}$ We use it to show that $T\in L(\ell_{2},X)$
given by $T(e_{n})=x_{n}, \, n\geq 1$, is 2-summing whenever
$(x_{n})\in \ell_{2}\{X\}.$ To this end, let $(\alpha^{(n)})$ be a
weakly 2-summable sequence in $\ell_{2}$. Then for $x^{*} \in
B_{X^{*}}$ and $n\geq 1$ , we see that
\begin{eqnarray}\nonumber
&&|<x^{*},T(\alpha^{(n)})>|= |<T^{*}(x^{*}),\alpha^{(n)}>| =
|{<x^{*},x_{i}>}_{i},\alpha^{(n)}|\\\nonumber&=&
|\sum_{i=1}^{\infty}<x^{*},x_{i}>\alpha_{i}^{(n)}| \leq
\sum_{i=1}^{\infty} |\alpha_{i}\alpha_{i}^{(n)}|  = c_{n}, say.
\end{eqnarray}

\noindent Finally, weak 2-summability of $(\alpha^{(n)}) $ in
$\ell_{2}$ yields that $c= (c_n)\in \ell_{2}$ and this completes the
argument. Regarding the proof of (ii), we note that  $(x_n)\in
\ell_{2}\{X_0\}$ if and only  if ST is a Hilbert-Schmidt map for
each $S\in L(X,H)$ where H is a Hilbert space, or equivalently,
${(ST)}^*=T^*S^*$ is a Hilbert-Schmidt map. Since S was chosen
arbitrarily, Remark 3.3(iv) yields that ${T}^*$ is 2-summing. \\

\noindent Now the identifications set up in (i) and (ii) above yield
that the equality of sequence spaces in (a) amounts to the
operator-ideal equality: $\Pi^{d_{2}}(\ell_{2},X)=
\Pi_{2}(\ell_{2},X) $  . However, it is well-known that the latter
equation holds precisely when $ X$ is Hilbertian (See \cite{9},
Theorem 4.19). An alternative argument based on the ``Eigenvalue
Theorem'' of Johnson et al \cite{15} proceeds as follows. Let $T:
X\rightarrow X$ be a nuclear map on $X$. Then $T$ can be factored as
$T = T_{2}D_{2}D_{1}T_{1}$ where $T_{1}\in
L(X,\ell_{\infty}),T_{2}\in L(\ell_{1},X), D_{1} \in
\Pi_{2}(\ell_{\infty},\ell_{2}), D_{2} \in
\Pi^{d_{2}}(\ell_{2},\ell_{1}) $. The ideal property of
$\Pi^{d}_{2}$ combined with the above equation gives $T_{2}D_{2} \in
\Pi_{2}(\ell_{2},X)$. Thus as a composite of two 2-summing maps, T
has absolutely summable eigenvalues by \cite{22}, Proposition$
3.4.5$ and Theorem $3.7.1$, and so X is isomorphic to a Hilbert
space, by the 'Eigenvalue Theorem'. See \cite{31} for further
applications of the Eigenvalue Theorem in the context of
factorization of certain operator ideal equations.\\

\noindent {\bf Proof (b):} It is easily seen that the equality
$\ell_{2} \{X_{0}\}= \ell_{2}[X]$ involving
 sequence spaces translates into the equality $L(X,\ell_{2})= \Pi_{2}(X,\ell_{2})$ involving operator
  ideals which, by virtue of Theorem $III.c.3$, holds exactly when X is a Hilbert Schmidt space.\\

\noindent (ii). The (FD)- property in question says that there exist
positive definite functions on each infinite
 dimensional Banach space which do not arise as the Fourier transform of a regular Borel measure on the
  dual of $X$. However, it turns out that if $X$ is a Hilbert space, there always exists a locally convex
   topology ù on $X$ ( the so-called Sazonov topology) such that each positive definite function on $X$ which
    is continuous in this topology is already a Fourier transform. Remarkably, it turns out that as soon as such
    a topology exists on a Banach space $X$, then $X$ is a Hilbert space. Thus, the locally convex topology $\tau$
     would provide the missing link if the following statement were true:\\
$(*)$ Every positive definite function on X is ù-continuous if and only if X is a Hilbert-Schmidt space.\\
However, it is not known if $(*)$ is true.\\

\noindent (iii). In the theory of vector measures, the so-called
‘localisation problem’ deals with the issue of enclosing
sequences from certain distinguished sequence spaces $\lambda (X)$
from a Banach space X inside ranges of vector measures( with or
without bounded variation) taking values in $X$. For $\lambda
(X)=c_{0}(X)$, it turns out that Banach spaces $ X$ for which
$\lambda (X)\subset R(X)$ are precisely those for which $X^*$ is a
subspace of an $L_{1}$− space. Here the missing object in the
desired factorization is a sequence space modeled on $X$ which lies
between $R(X)$ and $R_{bv}(X)$ and yields properties $(Q)$ and $
(R)$ for which $(P) = (Q) \wedge (R)$. We define
$$R_{vbv}(X)=\{ (x_{n})\subset X: \exists\,\, Banch\, space\,\, Z\supset X\, and \, \mu \in M_{bv}(Z)\,
 s.\,t. (x_{n})\subset\,\,\, rg(\mu)\}$$
It can be proved that $R_{bv}(X)\subset R_{vbv}(X)\subset R(X)$. Now
the inclusion relation in (iii) is equivalent to: (a) $c_0(X)\subset
R_{vbv}(X)$ and (b) $R_{vbv}(X)\subset R_{bv}(X)$. It was proved by
Pineiro \cite{23}(see also \cite{27} for an alternative proof) that
(a) holds exactly if the underlying space is Hilbertian whereas (b)
holds if and only if N∗ is a (GT)- space (i.e. satisfies
Grothendieck's theorem: $L(X,\ell_{2})=\Pi_1(X,\ell_{2})$, a
property which is obviously stronger than that of being a
Hilbert-Schmidt space). In particular, $X$ is a Hilbert-Schmidt
space. A far reaching refinement of (iii) was proved by the author
\cite{29}, replacing the space $c_0(X)$ by the smaller sequence
space $\ell_p\{X\} \,for\, p>2$\\

\noindent (iv). It is well known that an operator acting between
Hilbert spaces is Hilbert-Schmidt if and only if its adjoint is. A
suitable analogue of this result in the Banach space setting shall
make sense once it is clear what the Banach space analogue of a
Hilbert-Schmidt map ought to look like. We have already noted (see
Remark 3.3 (iii)) that p-summing maps coincide with Hilbert-Schmidt
maps on Hilbert spaces, with the equivalence of the p-summing norm
and the Hilbert-Schmidt norm. However, considering that the natural
norm on the ideal of 2-summing maps even coincides with the
Hilbert-Schmidt norm, the class of 2-summing maps stands out as the
most appropriate candidate for the Banach space analogue of
Hilbert-Schmidt maps. In view of this,
 it is natural to ask if the above stated result is valid for 2-summing maps in the Banach space
  setting. It turns out that finite dimensional spaces are the only Banach spaces X for which this
   property holds, i.e. such that (iv) holds. The desired factorization of (iv), therefore, amounts to the inclusions:\\

\noindent (a) $\Pi_{2}(X,Z)\subset \Pi^{d_{2}}(X,Z)\,\, \forall\,\, Banach\,spaces Z. $\\
(b)  $\Pi^{d_{2}}(X,Z)\subset \Pi_{2}(X,Z)\,\,\forall\,\, Banach\,spaces Z. $\\

\noindent As has been seen to be the case in respect of the
(FD)-properties encountered earlier, the above inclusions are valid
precisely when the Banach space X is Hilbert and Hilbert-Schmidt,
respectively. The proof makes use of Grothendieck's theorem quoted
earlier combined with some important estimates involving p-summing
maps. For details, see \cite{9}, Chapter 4.\\

\noindent (v).The fact that this condition on X is an (FD)-property
was proved by the author in \cite{27} where similar other results
are proved in connection with some problems arising in the theory of
vector measures. The missing object in the above factorization
entails the search for an ideal $Α$ of operators which can be
'sandwiched' between $\Pi_{2}$ and $N$ such that the resulting
equations characterize $(Q)$ and $(R)$,
 respectively. More precisely, we shall 'invent' an operator ideal $A$ such that\\
(i) $\Pi_{2}(X,\ell_{2}) = A(X,\ell_{2})$ iff X is Hilbertian.\\
(ii)$ A(X,\ell_{2})= N(X,\ell_{2} )$iff X is Hilbert-Schmidt.\\

\noindent Unfortunately, the natural choice for A, namely the ideal
$\Pi_{1}$ of absolutely summing maps that suggests itself for
effecting the desired factorization does not work$!$ In fact, it was
proved in \cite{27} that the class of Banach spaces  $X $ satisfying
the equation $\Pi_{2}(X,\ell_{2})= N(X,\ell_{2} )$
 are precisely those for which $X^*$ has the Gordon-Lewis
property and satisfies the Grothendieck theorem. However, It can be
shown (see \cite{31}) that we can choose the 'missing link' to be $A
= \Pi_{HS}$ in order to achieve the desired factorization of the
(FD)-property (P) given in Theorem(8) above as $(P) = (Q) \wedge
(R)$ where the properties  $(Q) $ and  $(R) $ are defined by the
above equations (i) and (ii). The operator ideal  $\Pi_{HS} $ is
defined by:$\Pi_{HS}(X,Y)=$
$$\{T\in L(X,Y); \exists\,T_{1}\in L(X,\ell_{2}), S \in  HS(\ell_{2},\ell_{2}),T_{2}
\in L(\ell_{2},Y)\,\, s.\, t.\, T= T_{2}ST_{1}\}$$

\begin{remark}
 By modifying the technique employed in the above factorization, it
is possible to factorise the (FD)-properties involving the following
operator-ideal
equations.\\
(a) $K(\ell_{2},X) = N(\ell_{2},X),$\\
(b) $S_{2}^{(e)}(\ell_{2},X)= N(\ell_{2},X).$\\

\noindent Here, $K$ denotes the class of compact maps whereas the
symbol $S_{2}^{(e)}(\ell_{2},X)$ is used for those maps for which
$\{T(x_{n})\}$ (absolutely) 2-summable
 in $X$ for some weakly 2-summable sequence $\{x_{n}\}$ in $\ell_{2}$. In each case, it turns out that
 the ideal $\Pi_{2}$ of 2-summing maps
provides the desired missing link! (See \cite{26}).
\end{remark}
\begin{theorem}
 For a Banach space, the following statements each in (A) and (B) are equivalent:\\
(A)(i) $X$ is a Hilbert space.\\
\indent (ii) $K(\ell_{2},X)=\Pi_{2}(\ell_{2},X).$\\
\indent (iii)$S_{2}^{(e)}(\ell_{2},X)=\Pi_{2}(\ell_{2},X).$\\
(B) (i) X is a Hilbert-Schmidt space.\\
\indent (ii)$\Pi_{2}(\ell_{2},X)=N(\ell_{2},X).$
\end{theorem}

\noindent {\bf Concluding Remarks:} The three important features of
(FD)-properties discussed as the main theme of this paper have been
illustrated with a number of examples drawn from different
situations witnessing finite dimensionality in functional  analysis.
However, notwithstanding the preponderance of examples cropping up
in our discussion of these phenomena, it should be emphasized that
not all (FD)-properties that one encounters in analysis lend
themselves to a suitable analogue under each of these categories.
For example, the set of objects failing the (FD)- property: $K(X) =
L(X)$ for an infinite dimensional Banach space X is not necessarily
'big', at least in certain pathological situations. In fact, it
turns out that for the Argyros-Haydon example of a Banach space
admitting a ‘small’ space of operators, the set $L(X)/K(X)$ does
not even contain a 2-dimensional space! In a similar vein, the
Heine-Borel property or the property: $X^*=X$΄does not characterize
nuclearity of the Frechet space  $X $. This motivates the search for
a set of theorems identifying a given (FD)-property as being
amenable to yield itself to one of the several features of finite
dimensionality as spelled out in the previous sections. The search
for such theorems promises to be a fruitful line of research in this
circle of ideas.

\noindent {\bf Acknowledgements:} This article is based on a series
of lectures the author had given at the Department of Mathematics,
I.I.T.Kanpur in Feb.2011. He would like to thank Prof. Manjul Gupta
of the department for arranging this visit and the Head of the
Department for his kind invitation and support to cover this visit
as a visiting faculty. It is also a pleasure to thank the anonymous
referees for their thorough analysis and useful remarks which have
helped improve the text.

\end{document}